\documentclass[10pt]{article}
\usepackage[cp1251]{inputenc}
\usepackage[english,russian]{babel}     %       Babel
\usepackage{amsfonts}
\usepackage{amsthm} % \newtheorem*
\usepackage{amssymb}
\usepackage{latexsym} % дополнительные символы

\newtheorem{teo}{\quad Theorem}

\newtheorem{df}{\quad Definition}

\newtheorem{vspteo}{\quad Theorem}
\newtheorem{vsplem}{\quad Lemma}
\renewcommand{\refname}{Literature}

\begin{document}
\Large

\noindent UDK 517.5

\noindent
{\bf S.O.~Chaichenko}

\noindent Donbass State Pedagogical University

\bigskip\noindent
{\bf APPROXIMATIONS OF PERIODIC FUNCTIONS BY \\ ANALOGUE OF ZIGMUND'S SUMS IN THE SPACES $L^{p(\cdot)}$}

\bigskip

{\noindent \footnotesize In this work we found order estimates for the upper bounds of the deviations of analogue of Zigmund's sums on the classes of $(\psi;\beta)$-differentiable functions in the metrics of generalized Lebesgue spaces with variable exponent.}

\bigskip

{\bf 1. Definition and formulation of the problem.} Let ${p}=p(x)$ be a $2\pi$-periodic measurable and essentially bounded function and let  $L^{p(\cdot)}$ be space of measurable
$2\pi$-periodic functions $f$ such that
$$
\int\limits_{-\pi}^\pi |f(x)|^{p(x)}~dx < \infty.
$$

If  $\underline{p}:={\rm ess}\inf\limits_{x} |p(x)|>1$ and
$\bar{p}:= {\rm ess}\sup\limits_{x} |p(x)|<\infty,$ then
$L^{{p(\cdot)}}$ are Banach spaces \cite{Sharapud.1979.Topolog}
(see., also \cite{Kovacik-Rakosnik}) with the norm, which can be
given by the formula
$$
    \|f\|_{p(\cdot)}:=\inf
    \Bigg\{{\alpha >0}:~ \int\limits_{-\pi}^\pi
    \Bigg|\frac{f(x)}{\alpha}  \Biggr|^{p(x)}~dx \le 1 \Bigg\}.
$$

Here are some definitions which will
used in the statement and proof of the results of this
article.

\begin{df}
It is said that a function $p=p(x)$ satisfies the Dini-Lipschitz condition of order
$\gamma,$ if
$$
    \omega (p; \delta) \Bigg(\ln \frac{1}{\delta} \Bigg)^\gamma \le K, \quad
    0<\delta<1,
$$
where
$$
    \omega (p; \delta)=\sup_{x_1, x_2 \in [-\pi;\pi]}
    \bigg\{|p(x_1)-p(x_2)|: \quad |x_1-x_2|\le \delta \bigg\}.
$$
\end{df}

The set of $2\pi$-periodic exponents $p=p(x)>1,$
satisfying  the Dini-Lipschitz condition of order $\gamma\ge 1$ in the period, is denoted by ${\cal P}^\gamma.$ Obviously, if $p \in {\cal P}^\gamma,$ then
$\underline{p}>1$ and $\bar{p}<\infty.$

In the work \cite{Sharapud.1979.Topolog} shown that when $1< \underline{p},~ \bar{p} <\infty,$ space
$L^{ q(\cdot)},$ where $q(x)=\frac{p(x)}{p(x)-1},$ is conjugate for
$L^{{p(\cdot)}}$ and for arbitrary functions $f \in L^{{p(\cdot)}}$ and $g\in
L^{ q(\cdot)}$  an analogue of the classical H\"{o}lder's inequality is true:
\begin{equation}\label{1}
    \int\limits_{-\pi}^\pi |f(x)g(x)|~dx \le K_{p,q}
    \|f\|_{{p(\cdot)}}\|g\|_{q(\cdot)}, \quad (K_{p,q}\le {1}/{\underline{p}}+{1}/{\underline{q}}),
\end{equation}
which, in particular, implies embedding: $L^{p(\cdot)} \subset L,$ where  $L$ is space of $2\pi$-periodic integrable on the period functions.

The space $L^{p(\cdot)}$ are called generalized Lebesgue spaces
with variable exponent. It is clear, that if $p=p(x)={\rm
const}>0,$ spaces $L^{{p(\cdot)}}$ coincide with the classical
Lebesgue spaces $L_p.$ In its turn, if $ \bar{p} <\infty,$
spaces $L^{{p(\cdot)}}$ are a special case of the so-called space
Orlicz-Musielak \cite{Musielak1983}.
For the first time, Lebesgue space with  variable exponent  appeared in the literature in the article W.~Orlicz \cite{Orlicz}. In the work \cite{Nakano} spaces  $L^{{p(\cdot)}}$
considered as an example of the more general function spaces and,
furthermore, been studied by many authors in different directions. The basic results of the theory of these spaces are available,
for example, in \cite{Sharapud.1979.Topolog,  Kovacik-Rakosnik,
Samko-1994, Sharap-1996-ravn-ogran, Fan-Zhao,
Sharap-2007-vopr-thory-aprox, Bandaliev_2014}. Note also that the generalized
Lebesgue spaces with variable exponent used in the theory
elastic mechanics, the theory of differential operators, variations
calculus \cite{Diening-Ruzicka, Ruzicka, Samko-2005}.

Next, we need the definitions of the $(\psi; {\beta})$-derivative and the sets
$L_{{\beta}}^\psi$, which belongs to A.I.~Stepanetz \cite[c.~142 --
143]{Stepanetz-2002m1}.

\begin{df}
Let $f\in L$ and
\begin{equation}\label{2}
    S[f]=\frac{a_0(f)}{2}+\sum_{k=1}^\infty(a_k(f)\cos~kx+b_k(f)\sin~kx)
    \equiv \sum_{k=0}^\infty A_k(f,x)
\end{equation}
be its Fourier series. Let, further, $\psi(k)$
be arbitrary function of natural argument and $\beta \in \mathbb{R}$.
Assume that the series
\begin{equation}\label{3}
    \sum_{k=1}^{\infty}\frac{1}{\psi(k)}
    \left( a_k(f)\cos\bigg(kx+\frac{\beta \pi}{2}\bigg)+
    b_k(f)\sin\bigg(kx+\frac{\beta \pi}{2}\bigg) \right)
\end{equation}
is the Fourier series of some function from $ L.$ This function is denoted by
$f_{{\beta}}^\psi(\cdot)$ (or $(D_{{\beta}}^\psi f)(\cdot)$) and
called $(\psi;{\beta})$-derivative of a function $f(\cdot).$ The set of
functions $f(\cdot),$ satisfying this condition is denoted by
$L_{{\beta}}^\psi.$
\end{df}

Denote by $L_{\beta,p(\cdot)}^{\psi}$
the classes $(\psi; {\beta})$-differentiable functions $f \in
L,$ such that $f_{{\beta}}^\psi \in L^{p(\cdot)}$, and by $\hat{Z}_n(f;x)$ the trigonometric polynomials of the form
\begin{equation}\label{4}
    \hat{Z}_n(f;x):=\frac{a_0(f)}{2}+\sum\limits_{k=1}^\infty
    \bigg(1-\frac{\psi(n)}{\psi(k)} \bigg) A_k(f;x).
\end{equation}

In this paper, we study the value
$$
    {\cal E} (L_{\beta,p(\cdot)}^{\psi};\hat{Z}_n)_{s(\cdot)}:=
    \sup_{f \in L_{\beta,p(\cdot)}^{\psi}} \|f- \hat{Z}_n(f)\|_{s(\cdot)}
$$
upper bounds of deviations analogues of Zigmund's sums $\hat{Z}_n(f;x)$ on the clas\-ses $L_{\beta,p(\cdot)}^{\psi}:=\{f \in L_{\beta}^{\psi}L^{p(\cdot)}: f^\psi_\beta \in U_{p(\cdot)} \}$, where $ U_{p(\cdot)}:=\{\varphi \in L^{p(\cdot)}: ~ \|\varphi\|_{p(\cdot)}\le 1\} $ --- the unit ball of $L^{p(\cdot)}.$

Note that in the case $\psi(k)=1/k^{r},~r>0,$ the sums (\ref{4}) are well known  Zigmund's sums
$$
        Z_n(f;x):=\frac{a_0(f)}{2}+\sum\limits_{k=1}^\infty
    \bigg(1-\bigg(\frac{k}{n}\bigg)^r \bigg) A_k(f;x).
$$

{\bf 2. Auxiliary results.}
In the proof of the main assertions of this work we use the following well-known results.

\begin{vspteo}\label{T.A}\cite{Sharap-2007-vopr-thory-aprox}
If $p\in {\cal P}^\gamma,$ then for an arbitrary function
$f \in L^{{p(\cdot)}}$ the inequalities hold
\begin{equation}\label{5}
    \| S_n(f)\|_{{p(\cdot)}} \le C_{p} \| f\|_{{p(\cdot)}},
\end{equation}
\begin{equation}\label{6}
    \| \tilde{f}\|_{{p(\cdot)}} \le K_{p} \| f\|_{{p(\cdot)}},
\end{equation}
where
$\tilde{f}(\cdot)$ is a functions, trigonometric conjugate to $f(\cdot)$, and  $C_{p}$ and $K_p$ are  positive constants which does not depend on $n$ and $f.$
\end{vspteo}

From the inequality (\ref{5}) of  Theorem \ref{T.A}, in particular,
follows that for an arbitrary function $f \in L^{{p(\cdot)}},$ on condition
$p\in {\cal P}^\gamma,$ its Fourier series converges to $f$ in the metric of the spaces $L^{{p(\cdot)}},$ that is
\begin{equation}\label{7}
    \| f-S_n(f) \|_{{p(\cdot)}} \to 0, \quad n \to \infty,
\end{equation}
and the relation holds
\begin{equation}\label{8}
    E_{n}(f)_{{p(\cdot)}} \le \| f-S_{n-1}(f) \|_{{p(\cdot)}} \le K_{p} E_{n}(f)_{{p(\cdot)}},
\end{equation}
where
$$
    E_{n} (\varphi)_{{p(\cdot)}}:=\inf\limits_{t_{n-1} \in
    {\cal T}_{n-1}} \| \varphi- t_{n-1}\|_{{p(\cdot)}},
    \quad \varphi \in L^{{p(\cdot)}},
$$
is the best approximation of $\varphi$ by subspace ${\cal T}_{n-1}$
trigonometric polynomials of order, not higher than $n-1,$ and $K_{p}$
is value which depends only on $p=p(x).$

\begin{vsplem}\label{L.A} \cite{Kokilash_Samko_2009}
Let the sequence $\mu(k),~k=0,1,2,\ldots,$ satisfies the conditions
$$
    \nu_0=\nu_0(\mu)=\sup_k|\mu(k)| \le C, \qquad
    \sigma_0=\sigma_0(\mu)=\sup_{m \in \mathbb{N}} \sum_{k=2^m}^{2^{m+1}} |\mu(k+1)-\mu(k)|\le C,
$$
where $C$ is value which does not depend on $k$ and $m.$

Then, if $p\in {\cal P}^\gamma,$ for a given function
$f \in L^{{p(\cdot)}}$ there exists a function $F\in L^{{p(\cdot)}}$  such that the series
$$
    \frac{\mu(0)a_0(f)}{2}+\sum_{k=1}^\infty\mu(k)(a_k(f)\cos~kx+b_k(f)\sin~kx)
$$
is the Fourier series of $F$ and the inequality is true
\begin{equation}\label{9}
    \|F\|_{{p(\cdot)}}\le K \lambda \|f\|_{{p(\cdot)}}, \quad \lambda=\max\{\nu_0, \sigma_0\},
\end{equation}
where the value $K$  does not depend on the function $f.$
\end{vsplem}

In the case $p=p(x)\equiv {\rm const}$ this statement is a well-known lemma of Marcinkievicz for multipliers \cite{Marcinkievicz}.

We will also use the following theorem of Hardy-Littlewood.

\begin{vspteo}\label{T.B}\cite{Hardy_Littlewood_1928}
Let $1<p<s<\infty,~p,s={\rm const},$ $\alpha=p^{-1}-s^{-1}$ and
$$
    D_\alpha(t)=\sum_{k=1}^\infty k^{-\alpha} \cos kt.
$$

Then, for an arbitrary function $\varphi \in L_p$ the convolution
$$
    \Phi_\alpha(x)=\frac{1}{\pi} \int\limits_{-\pi}^\pi \varphi (x+t) D_\alpha(t)~dt
$$
belongs to $L_s,$ and
$$
\|\Phi_\alpha\|_s \le C_{s,p} \| \varphi \|_p,
$$
where the value $C_{s,p}$  depends only on $s$ and $p.$
\end{vspteo}

Note that if $\varphi \in L_p$ and  $S[\varphi]=\sum\limits_{k=0}^\infty A_k(\varphi;x),$ then
$$
    S[\Phi_\alpha]=\sum\limits_{k=0}^\infty k^{-\alpha} A_k(\varphi;x),
$$
that is $\Phi_\alpha=M_\alpha (\varphi),$ where $M_\alpha$ --- operator-multiplier, which is determined by the sequence of $\mu_\alpha(k)=k^{-\alpha},~k=0,1,2,\ldots,$ and it acts from $L_p$ to $L_s,$ where indicators $1<p<s<\infty,~p,s={\rm const},$  are related by the equation $p^{-1}-s^{-1}=\alpha.$

{\bf 2. Approximation by analogue of Zigmund's sums.}

We define the function $\mu_{n,\alpha}(k)$ and $\tilde{\mu}_{n,\alpha}(k)$ as follows:
\begin{equation}\label{10}
    \mu_{n,\alpha}(k)=\cases{k^\alpha \psi(n) \cos \frac{\beta \pi}{2}, & $1 \le k \le n-1,$ \cr
                            k^\alpha \psi(k) \cos \frac{\beta \pi}{2}, & $n \le k;$}
\end{equation}
\begin{equation}\label{11}
    \tilde{\mu}_{n,\alpha}(k)=\cases{k^\alpha \psi(n) \sin \frac{\beta \pi}{2}, & $1 \le k \le n-1,$ \cr
                            k^\alpha \psi(k) \sin \frac{\beta \pi}{2}, & $n \le k;$}
\end{equation}

For each fixed $\alpha \ge 0$ we denote by $\Upsilon_{\alpha, n}$  the set of pairs $(\psi;\beta)$, such that for any positive number $n$  the conditions
\begin{equation}\label{12}
    \nu_\alpha(\psi; \beta; n):=\sup_{k} |\mu_{n,\alpha}(k)| \le C \nu(n)n^\alpha<K,
\end{equation}
\begin{equation}\label{13}
    \sigma_\alpha (\psi;\beta; n) := \sup_{m \in \mathbb{N}} \sum_{k=2^m}^{2^{m+1}}
    |\mu_{n,\alpha}(k+1)-\mu_{n,\alpha}(k)| \le C \nu (n) n^\alpha< K,
\end{equation}
and similar conditions for the function $\tilde{\mu}_{n,\alpha}(k)$ hold, where $\nu(n)=\sup_{k\ge n} |\psi(k)|,$  $C$ and $K$ are positive constants uniformly bounded on $n.$

At first we consider  the case when the functions $p=p(x)$ and $s=s(x)$ on the period satisfy  the inequality $ s(x)\le p(x).$ In our notation, the following assertion is true.

\begin{teo}\label{T.1}
Let $(\psi;\beta) \in \Upsilon_{0,n}$ и $ p, s \in \mathcal{P}^\gamma,$  $ s(x)\le p(x),~x\in [0;2\pi].$ Then, for all $n \in \mathbb{N}$ the inequality
$$
    C_{p,s}\nu(n) \le {\cal E} (L_{\beta, p(\cdot)}^\psi; \hat{Z}_n)_{ s(\cdot)} \le K_{p,s}\nu (n),
$$
holds, where $C_{p,s}$ and $K_{p,s}$ are some constants that do not depend on  $n$.

\end{teo}

{\bf Proof.} Suppose first that $p(x)\equiv s(x),~x\in [0;2\pi].$ For an arbitrary function $f \in L_{\beta,p(\cdot)}^{\psi},$ the equality is true
$$
    f(x)-\hat{Z}_n(f;x)=\sum_{k=1}^{n-1} \frac{\psi(n)}{\psi(k)} A_k(f;x)+
    \sum_{k=n}^\infty  A_k(f;x)=
$$
\begin{equation}\label{14}
    =\sum_{k=1}^\infty \mu_{n,0} A_k(f^\psi_\beta; x)+
    \sum_{k=1}^\infty \tilde{\mu}_{n,0} \tilde{A}_k(f^\psi_\beta; x)
    := M_0(f^\psi_\beta)+\tilde{M}_0(\tilde{f}^\psi_\beta),
\end{equation}
where $M_0$ and $\tilde{M}_0$ are operators-multipliers, which are defined by the sequences (\ref{10}) and (\ref{11}) respectively, $\alpha=0.$

According to the conditions of the theorem, the couples $(\psi;\beta)$ belong to the set $\Upsilon_{0, n},$ therefore, the sequence (\ref{10}) and (\ref{11}) satisfy the conditions of lemma \ref{L.A}. Applying this lemma, given the inequalities (\ref{6}), (\ref{12}) and (\ref{13}), for an arbitrary function $f \in L_{\beta,p(\cdot)}^{\psi}$ on the basis of the equality (\ref{14}) we find
$$
    \|f(\cdot)-\hat{Z}_n(f;\cdot) \|_{p(\cdot)} =
    \|M_0(f^\psi_\beta)+\tilde{M}_0(\tilde{f}^\psi_\beta) \|_{p(\cdot)}\le
$$
\begin{equation}\label{15}
    \le K\nu(n)(\| f^\psi_\beta\|_{p(\cdot)} + \| \tilde{f}^\psi_\beta\|_{p(\cdot)})\le C_p \nu(n),
\end{equation}
where $C_p$ is positive constant that depends only on the function $p=p(x).$

In the article \cite{Sharapudinov_O_bazist_sist_Haara} was shown that if $1\le s(x)\le p(x) \le \bar{p}<\infty,$ then for an arbitrary function $f \in L^{p(\cdot)}$ the inequality hold
\begin{equation}\label{16}
    \|f\|_{s(\cdot)}\le K_{s,p}\|f\|_{p(\cdot)}.
\end{equation}

From the relations (\ref{15}) and (\ref{16}) we obtain the estimate
\begin{equation}\label{17}
    {\cal E} (L_{\beta,p(\cdot)}^{\psi};\hat{Z}_n)_{s(\cdot)} \le
    {\cal E} (L_{\beta,s(\cdot)}^{\psi};\hat{Z}_n)_{s(\cdot)} \le C_{p,s} \nu (n),
\end{equation}
where $C_{p,s}$ is positive constant that depends only on the functions $p=p(x)$ and $s=s(x).$

We now obtain the lower estimate for the value of ${\cal E} (L_{\beta,p(\cdot)}^{\psi};\hat{Z}_n)_{s(\cdot)}.$ If for given function $\psi(k)$ and the number $n\in \mathbb{N}$ there exist the natural number $k_n,$ for which the equality
\begin{equation}\label{18}
    \nu (n)=\sup_{k\ge n} | \psi(k)|=\psi(k_n),
\end{equation}
is true, then the corresponding lower estimate can be obtained with help of the function
\begin{equation}\label{19}
    f_n(x)=\frac{\psi(k_n)}{\|\cos k_n x\|_{p(\cdot)}} \cos (k_n x-\frac{\beta\pi}{2}).
\end{equation}

Indeed, since
$$
    \| (f_n(x))^\psi_\beta\|_{p(\cdot)}= \bigg\| \frac{\cos k_n x}
    {\|\cos k_n x\|_{p(\cdot)}} \bigg\|_{p(\cdot)}=1,
$$
then the function $f_n \in L_{\beta,p(\cdot)}^{\psi}$ and
$$
    {\cal E} (L_{\beta,p(\cdot)}^{\psi};\hat{Z}_n)_{s(\cdot)}\ge \|f_n-\hat{Z}(f_n) \|_{s(\cdot)}=
$$
$$
    =\frac{\psi(k_n)}{\|\cos k_n x\|_{p(\cdot)}} \| \cos k_n x\|_{s(\cdot)}=C_{p,s} \nu(n).
$$

But if for the function $\psi(k)$ and the number $n\in \mathbb{N}$
there not exist a natural number $k_n,$ for which the equality (\ref{18}) holds,  due to the limitation of the set  $\{|\psi(k)|\}$ of values of the function $\psi(k)$ we will have
$$
    \nu (n)=\sup_{k\ge n} | \psi(k)| = \sup_{k\ge n} \{|\psi(k)|\}.
$$
In this case, there exists a sequence $k_j,~j\in \mathbb{N}$ such that $k_j\ge n$ and the numbers  $\psi(k_j)$ don't decrease and converge to $\nu(n).$ Let $\mathbb{A}=\cup_{j} f_j(x),$
where the function $f_j(x)$
defined by the equality (\ref{19}). Since $f_j \in L_{\beta,p(\cdot)}^{\psi}$ for any $j\in \mathbb{N},$ than
$$
    {\cal E} (L_{\beta,p(\cdot)}^{\psi};\hat{Z}_n)_{s(\cdot)}=
    \sup_{f \in L_{\beta,p(\cdot)}^{\psi}} \|f-\hat{Z}_n(f) \|_{s(\cdot)}\ge
    \sup_{f \in \mathbb{A}} \|f-\hat{Z}_n(f) \|_{s(\cdot)}=
$$
$$
     =\sup_{j \in \mathbb{N}} \frac{\psi(k_j)}{\|\cos k_j x\|_{p(\cdot)}}
     \| \cos k_j x\|_{s(\cdot)}=C_{p,s} \nu(n).
$$

The theorem is proved.

We now obtain an estimate of the sequence ${\cal E} (L_{\beta,p(\cdot)}^{\psi};\hat{Z}_n)_{s(\cdot)}$
in the case where the function $p=p(x)$ and $s=s(x)$ on the period satisfy the inequality $p(x)<s(x).$  The following result gives the upper estimate.

\begin{teo}\label{T.2}
Let $ p, s \in \mathcal{P}^\gamma,$  $ p(x)< s(x),~x\in [0;2\pi]$ and $(\psi;\beta) \in \Upsilon_{\alpha,n},$ $\alpha={1}/{\underline{p}}-{1}/{\overline {s}}.$
Then, for all $n \in \mathbb{N}$  the following inequality
\begin{equation}\label{20}
   {\cal E} (L_{\beta, p(\cdot)}^\psi; \hat{Z}_n)_{ s(\cdot)} \le C_{p,s} n^\alpha \nu (n),
\end{equation}
hold, where  $C_{p,s}$  is a positive constant which independent of $n$.

\end{teo}

{\bf Proof.}
For an arbitrary function $f \in L_{\beta,p(\cdot)}^{\psi},$ the equality
$$
    f(x)-\hat{Z}_n(f;x)=\sum_{k=1}^{n-1} \frac{\psi(n)}{\psi(k)} A_k(f;x)+
    \sum_{k=n}^\infty  A_k(f;x)=
$$
\begin{equation}\label{21}
    =\sum_{k=1}^\infty \mu_{n,\alpha} k^{-\alpha} A_k(f^\psi_\beta; x)+
    \sum_{k=1}^\infty \tilde{\mu}_{n,\alpha} k^{-\alpha}\tilde{A}_k(f^\psi_\beta; x)
    := M_\alpha (g_\alpha)+\tilde{M}_\alpha(\tilde{g}_\alpha),
\end{equation}
holds, where $M_\alpha$ and $\tilde{M}_\alpha$ are operators-multipliers, which are defined by the sequences (\ref{10}) and (\ref{11}) respectively, $\alpha={1}/{\underline{p}}-{1}/{\overline {s}}$ and
$$
    g_\alpha(x):=\sum_{k=1}^\infty k^{-\alpha} A_k(f^\psi_\beta; x)=
    \frac{1}{\pi} \int\limits_0^{2\pi} f^\psi_\beta(x+t) D_\alpha(t)~dt,
$$
$$
  \tilde{g}_\alpha(x):=\sum_{k=1}^\infty k^{-\alpha} \tilde{A}_k(f^\psi_\beta; x)=
    \frac{1}{\pi} \int\limits_0^{2\pi} \tilde{f}^\psi_\beta(x+t) D_\alpha(t)~dt,
$$
$D_\alpha(t)$ is function defined in Theorem \ref{T.B}.

Since $f \in L^{\psi}_{\beta,p(\cdot)},$ then $f^\psi_\beta \in L^{p(\cdot)}$, and
moreover  $f^\psi_\beta \in L^{\underline{p}}.$ By Theorem \ref{T.B}
we conclude that the convolution $g_\alpha(x)$  belongs $L^{\overline{s}},$ and moreover $g_\alpha \in L^{s(\cdot)}.$ From the condition $(\psi;\beta) \in \Upsilon_{n,\alpha}$ by lemma \ref{L.A} we conclude that the operator-multiplier $M_\alpha$ acts from $L^{s(\cdot)}$ to $L^{s(\cdot)}$ for any $s\in \mathcal{P}^\gamma.$ Using analogous arguments for the function $\tilde{g}_\alpha(x),$ taking into account  the inequalities (\ref{6}), (\ref{9}),  (\ref{12}) and (\ref{13}), for an arbitrary function $f \in L_{\beta,p(\cdot)}^{\psi}$ on the basis of the equality (\ref{21}) we find
$$
    \|f-\hat{Z}_n(f)\|_{s(\cdot)} \le \|M_\alpha (g_\alpha)\|_{s(\cdot)}+\|\tilde{M}_\alpha(\tilde{g}_\alpha)\|_{s(\cdot)}\le
   K n^\alpha \nu(n) (\|g_\alpha\|_{s(\cdot)}+\|\tilde{g}_\alpha\|_{s(\cdot)}) \le
$$
$$
    \le C_{p,s} n^\alpha \nu(n) ( \| f^\psi_\beta \|_{p(\cdot)}+ \| \tilde{f}^\psi_\beta\|_{p(\cdot)})\le
    C_{p,s} n^\alpha \nu(n).
$$

To make formulate the following assertion, which gives a lower estimate for the quantity
${\cal E} (L_{\beta,p(\cdot)}^{\psi};\hat{Z}_n)_{s(\cdot)}$ in the case, if the function $p=p(x)$ and $s=s(x)$  satisfy the inequality $p(x)<s(x)$ on the period, we need following definition.

Denoted by $\mathfrak{B}$ the set of pairs $(\psi;\beta),$ such that for any $n \in \mathbb{N}$ the relations are true:
$$
    \sup\limits_{n \le k \le 2n} \bigg| \frac{\nu(n)}{\psi(k)}\bigg| \le C, \quad
    \sup_{m\in \mathbb{N}} \sum_{k=2^m}^{2^{m+1}} |\tau (k+1)-\tau(k)| \le C,
$$
where $C$ is a positive constant which independent of $n,$  $\nu(n)=\sup_{k\ge n} \psi(k)$ and
\begin{equation}\label{25}
    \tau(k):= \cases{0,  &  $1\le k \le n-1,$ \cr
                    \frac{\nu(n)}{\psi(k)}, & $n\le k \le 2n.$}
\end{equation}

\begin{teo}\label{T.3}
Let $ p, s \in \mathcal{P}^\gamma,$  $ p(x)< s(x),~x\in [0;2\pi]$ and $(\psi;\beta) \in \mathfrak{B}.$
Then, for all $n \in \mathbb{N}$ the following inequality is true
\begin{equation}\label{22}
  {\cal E}_n (L_{\beta, p(\cdot)}^\psi; \hat{Z}_n)_{ s(\cdot)} \ge
  C_{p,s} \nu (n) n^{{1}/{\overline{p}}-{1}/{\underline {s}}},
\end{equation}
where $C_{p,s}$  is a positive constant which independent of $n$.

\end{teo}

{\bf Proof.} For obtaining a lower estimate, let us show that for any positive integer $n$  in class $L_{\beta, p(\cdot)}^\psi$ there exists a function $f_n^*,$ for which the inequality is true
$$
    \| f_n^* - \hat{Z}_n\|_{s(\cdot)} \ge C_{p,g}   \nu(n) n^{{1}/{\overline{p}}-{1}/{\underline {s}}}.
$$

For this we fix $n \in \mathbb{N}$ and consider the function
$$
    f_n^*(x)=\sum\limits_{k=n}^{2n} \psi(k)\cos (kx-\frac{\beta \pi}{2}).
$$
Since
$$
    (f_n^*(x))^\psi_\beta=\sum\limits_{k=n}^{2n} \cos kx=\frac{\sin n x/2 \cos (3n+1)x/2}{\sin x/2},
$$
then using the relation (\ref{16}) and also well-known inequality
\begin{equation}\label{23}
    \frac{x}{\pi} \le \sin x/2,\quad \sin x \le x, \quad x \in [0;\pi],
\end{equation}
we obtain
$$
    \| (f_n^*)^\psi_\beta\|_{p(\cdot)}=\bigg\| \sum\limits_{k=n}^{2n} \cos kx\bigg\|_{p(\cdot)}\le
    K_p\bigg\| \sum\limits_{k=n}^{2n} \cos kx\bigg\|_{\overline{p}}  =
$$
$$
    =\bigg(2 \int\limits_0^{\pi} \bigg|\sum\limits_{k=n}^{2n} \cos kx \bigg|^{\overline{p}}dx
    \bigg)^{1/\overline{p}}  \le \bigg(2 \int\limits_0^{\pi} \bigg|\frac{\sin n x/2}{\sin x/2}
    \bigg|^{\overline{p}} dx\bigg)^{1/\overline{p}} \le C^*_p n^{1-1/\overline{p}}.
$$

This implies that the function
$$
    g^*_n(x)=\frac{n^{1/\overline{p}-1}}{C^*_p} f_n^*(x)=
    \frac{n^{1/\overline{p}-1}}{C^*_p} \sum\limits_{k=n}^{2n} \psi(k)\cos (kx-\frac{\beta \pi}{2})
$$
belongs to the class  $L_{\beta, p(\cdot)}^\psi$.

Again using the inequalities (\ref{16}) and (\ref{23}), we find
$$
     \bigg\| \sum\limits_{k=n}^{2n} \cos kx\bigg\|_{s(\cdot)} \ge
      \bigg(2 \int\limits_0^{\pi} \bigg|\frac{\sin n x/2 \cos(3n+1)x/2}{\sin x/2}
    \bigg|^{\underline{s}}dx \bigg)^{1/\underline{s}} \ge
$$
\begin{equation}\label{24}
    \ge C_{s} n^{1-1/\underline{s}} \bigg(\int\limits_0^{\pi/2}
    (\cos x)^{\underline{s}}~dx \bigg)^{1/\underline{s}}\ge
    K_{s} n^{1-1/\underline{s}}.
\end{equation}

If now by $T_\psi$ we denote the operator-multiplier that generates a sequence (\ref{25}),
then by lemma \ref{L.A} on condition $(\psi;\beta) \in \mathfrak{B}$ we will have
$$
  \bigg\| \sum\limits_{k=n}^{2n} \cos (kx-\frac{\beta \pi}{2}) \bigg\|_{s(\cdot)}=
  \bigg\| T_\psi \bigg(\sum\limits_{k=n}^{2n} \frac{\psi(k)}{\nu(n)}
  \cos (kx-\frac{\beta \pi}{2})\bigg)\bigg\|_{s(\cdot)}\le
$$
$$
    \le C \bigg\|\sum\limits_{k=n}^{2n} \frac{\psi(k)}{\nu(n)}
  \cos (kx-\frac{\beta \pi}{2})\bigg\|_{s(\cdot)}.
$$

Hence, considering the inequality (\ref{24}) we find
$$
    \bigg\|\sum\limits_{k=n}^{2n} \frac{\psi(k)}{\nu(n)}
    \cos (kx-\frac{\beta \pi}{2})\bigg\|_{s(\cdot)} \ge
    K \bigg\| \sum\limits_{k=n}^{2n} \cos (kx-\frac{\beta \pi}{2}) \bigg\|_{s(\cdot)}\ge
$$
\begin{equation}\label{26}
    \ge C_s \bigg\| \sum\limits_{k=n}^{2n} \cos kx \bigg\|_{\underline{s}}
    \ge K_{s} n^{1-1/\underline{s}}.
\end{equation}

Using the relation (\ref{26}), we obtain
$$
    {\cal E}_n (L_{\beta, p(\cdot)}^\psi; \hat{Z}_n)_{ s(\cdot)} \ge
    \| g^*_n - \hat{Z}_n(g^*_n)\|_{s(\cdot)}= \bigg\|\frac{n^{1/\overline{p}-1}}{C^*_p}
    \sum\limits_{k=n}^{2n} \psi(k)\cos (kx-\frac{\beta \pi}{2}) \bigg\|_{s(\cdot)}\ge
$$
$$
    \ge \frac{n^{1/\overline{p}-1}}{C^*_p} \nu(n) \bigg\|
    \sum\limits_{k=n}^{2n} \frac{\psi(k)}{\nu(n)}\cos (kx-\frac{\beta \pi}{2}) \bigg\|_{s(\cdot)}\ge
$$
$$
    \ge C_{p,s} n^{1/\overline{p}-1} \nu(n) n^{1-1/\underline{s}}=C_{p,s} \nu(n)
    n^{{1}/{\overline{p}}-{1}/{\underline {s}}}.
$$

The theorem is proved.

\renewcommand{\refname}{References}

\end{document}